\documentclass{elsart3-1}

\usepackage{amssymb}
\usepackage[T1]{fontenc}
\usepackage[latin1]{inputenc}
\usepackage{latexsym}
\usepackage{xypic}
\usepackage{eufrak}
\usepackage{euscript}
\usepackage{amsmath}
\usepackage{verbatim}
\usepackage{fancyhdr}

\usepackage[english,francais]{babel}

\newtheorem{theoreme}{Th\'eor\`eme}[section]
\newtheorem{lemm}[theoreme]{Lemme}

\newtheorem{defi}[theoreme]{D\'efinition\rm}

\def\og{\leavevmode\raise.3ex\hbox{$\scriptscriptstyle\langle\!\langle$~}}
\def\fg{\leavevmode\raise.3ex\hbox{~$\!\scriptscriptstyle\,\rangle\!\rangle$}}

\def\Cu{\EuScript{C}^1}
\def\Cun{\EuScript{C}^n}

\renewcommand{\d}{\displaystyle}

\renewcommand{\leq}{\leqslant}
\newcommand{\K}{\mathbb{K}}
\newcommand{\N}{\mathbb{N}}

\newcommand{\R}{\mathbb{R}}

\newcommand{\C}{\mathbb{C}}

\def\cal{\mathcal}

\begin{document}
\centerline{Probability Theory/Dynamical Systems}
\begin{frontmatter}




%
\selectlanguage{francais}
\title{\textbf{Plongement stochastique des syst\`{e}mes lagrangiens}}



\author[cresson]{Jacky CRESSON},
\ead{cresson@math.univ-fcomte.fr}
\author[darses]{S\'{e}bastien DARSES}
\ead{darses@math.univ-fcomte.fr}

\address[cresson]{Laboratoire de Math\'{e}matiques, Universit\'{e} de Franche-Comt\'{e}}
\address[darses]{Laboratoire de Math\'{e}matiques, Universit\'{e} de Franche-Comt\'{e}}


\medskip
\selectlanguage{francais}

\begin{abstract}
\selectlanguage{francais} On d\'efinit un op\'erateur agissant sur
des processus stochastiques, qui \'{e}tend la d\'{e}rivation classique sur
les fonctions d\'{e}terministes diff\'{e}rentiables. On utilise cet
op\'{e}rateur pour d\'{e}finir une proc\'{e}dure associant aux op\'{e}rateurs
diff\'{e}rentiels et \'{e}quations diff\'{e}rentielles ordinaires leurs
analogues stochastiques. Elle est appel\'{e}e plongement stochastique.
En plongeant les syst\`{e}mes lagrangiens, nous obtenons une \'{e}quation
d'Euler-Lagrange stochastique, qui dans le cas des syst\`{e}mes
lagrangiens naturels est appel\'{e}e \'{e}quation de Newton plong\'{e}e. Cette
derni\`{e}re contient l'\'{e}quation de Newton stochastique introduite par
Nelson dans sa th\'{e}orie dynamique des diffusions browniennes.
Enfin, on consid\`{e}re une diffusion \`{a} drift gradient, \`{a}
coefficient de diffusion constant et poss\'{e}dant une densit\'{e} de
probabilit\'{e}. On d\'{e}montre alors qu'une condition n\'{e}cessaire
pour que cette diffusion soit solution de l'\'{e}quation de Newton
plong\'{e}e, est que sa densit\'{e} soit le carr\'{e} du module d'une
fonction d'onde solution d'une \'{e}quation de Schr\"{o}dinger
lin\'{e}aire.

{\it Pour citer cet article~: A. Nom1, A. Nom2, C. R.
Acad. Sci. Paris, Ser. I 340 (2005).}
\vskip 0.5\baselineskip

\selectlanguage{english} \noindent{\bf Abstract} \vskip
0.5\baselineskip \noindent {\bf Stochastic embedding of lagrangian
systems. } We define an operator which extends classical
differentiation from smooth deterministic functions to certain
stochastic processes. Based on this operator, we define a
procedure which associates a stochastic analog to standard
differential operators and ordinary differential equations. We
call this procedure stochastic embedding. By embedding lagrangian
systems, we obtain a stochastic Euler-Lagrange equation which, in
the case of natural lagrangian systems, is called the embedded
Newton equation. This equation contains the stochastic Newton
equation introduced by Nelson in his dynamical theory of brownian
diffusions. Finally, we consider a diffusion with a gradient
drift, a constant diffusion coefficient and having a probability
density function. We prove that a necessary condition for this
diffusion to solve the embedded Newton equation is that its
density be the square of the modulus of a wave function solution
of a linear Schr\"{o}dinger equation.

 {\it To cite this article: A. Nom1, A. Nom2, C. R.
Acad. Sci. Paris, Ser. I 340 (2005).}
\end{abstract}
\end{frontmatter}




\selectlanguage{francais}

\section{D\'eriv\'ee stochastique dynamique}
\label{}

On note $I:=]a,b[$ o\`u $a<b$ et $J:=[a,b]$ l'adh\'erence de $I$
dans $\R$. Soit $\K$ un corps et $d\in\N^*$. On se donne un espace
probabilis\'e $(\Omega,\cal{A},P)$ sur lequel existent une famille
croissante de tribus $(\cal{P}_t)_{t\in J}$ et une famille
d\'ecroissante de tribus $(\cal{F}_t)_{t\in J}$. Suivant Yasue
\cite{ya1}, on introduit la

\begin{defi}

On note $\Cu_{\K}(J)$ l'ensemble des processus $X$ d\'efinis sur
$J\times \Omega$, \`a valeurs dans $\K^d$ et tels que : $X$ soit
$(\cal{P}_t)$ et $(\cal{F}_t)$ adapt\'{e}, pour tout $t\in J$ $X_t\in
L^2(\Omega)$, l'application $t\to X_t$ de $J$ dans $L^2(\Omega)$
est continue, pour tout $t\in I$ les quantit\'es
$DX_t=\lim_{h\rightarrow 0^+}h^{-1} E[X_{t+h}-X_t\mid {\cal P}_t
]$, et $D_* X_t=\lim_{h\rightarrow 0^+} h^{-1} E[X_t-X_{t-h}\mid
{\cal F}_t ]$, existent dans $L^2(\Omega)$, et enfin les
applications $t\to DX_t$ et $t\to D_*X_t$ sont continues de
$I$ dans $L^2(\Omega)$.\\
Le compl\'et\'e de $\Cu_{\K} (J)$ pour la norme $\parallel
X\parallel=\sup_{t\in I} (\parallel X_t\parallel_{L^2(\Omega)}
+\parallel DX_t\parallel_{L^2(\Omega)} +\parallel D_*
X(t)\parallel_{L^2(\Omega)} )$, est encore not\'e $\Cu_{\K}(J)$,
et simplement $\Cu(J)$ quand $\K=\R$.
\end{defi}

Les quantit\'{e}s $D$ et $D_*$ sont introduites par Edward Nelson dans
sa th\'{e}orie dynamique des diffusions browniennes
(\textit{cf}\cite{ne1}). Soit $\iota$ l'injection $\d \iota :
\left\{\begin{array}{ccc}
C^1(J) & \longrightarrow & \Cu (J) \\
f & \longmapsto & \iota(f) : (\omega,t)\mapsto f(t)
\end{array}\right.$. On note : $\cal{P}_{det}^k:=\iota(C^k(J))$. Le probl\`{e}me d'extension
consiste \`{a} trouver un op\'erateur $\cal{D}: \Cu(I)\to \Cu_{\C}(I)$
satisfaisant :
\begin{itemize}
    \item[(i)] (Recollement) $\d \cal{D}\iota(f)_t=\frac{df}{dt}(t)$ sur $\Omega$,
    \item[(ii)] ($\R$-lin\'earit\'e) $\cal{D}$ est $\R$-lin\'eaire,
    \item[(iii)] (Reconstruction) Si l'on note ${\cal D}X =A(DX ,D_* X)+iB(DX,D_* X),$ o\`u $A$ et $B$ sont des $\R$-formes
lin\'eraires, on suppose que l'application de $\R^2$ dans $\R^2$ : $(x,y)
\mapsto (A(x,y),B(x,y))$ est inversible.
\end{itemize}

L'op\'{e}rateur $\cal{D}$ \'{e}tend donc la d\'{e}rivation classique
(\textit{cf}(i)) en un op\'{e}rateur lin\'{e}aire (\textit{cf}(ii)) sur
$\Cu(J)$ et la connaissance de $\cal{D}X$ induit celle de $DX$ et
$D_*X$ (\textit{cf}(iii)). On obtient :

\begin{lemm} Les seuls op\'{e}rateurs $\Cu(I)\to \Cu_{\C}(I)$ v\'{e}rifiant (i), (ii) et (iii)
sont: $${\cal D}_{\mu} =\d {D+D_* \over 2} +i\mu {D-D_* \over 2} ,\ \mu =\pm 1
.$$
\end{lemm}

On \'{e}crira $\cal{D}:=\cal{D}_1$ et
$\overline{\cal{D}}:=\cal{D}_{-1}$. La d\'{e}finition des it\'{e}r\'{e}s de
$\cal{D}$ et $\overline{\cal{D}}$ n\'{e}cessite l'extension de ces
op\'{e}rateurs aux processus \`{a} valeur complexe. Dans la suite, on
\'{e}tend $\cal{D}$ et $\overline{\cal{D}}$ par $\C-$lin\'{e}arit\'{e}
aux processus complexes, \textit{i.e.} pour tous $X,Y\in\Cu(J)$,
$\cal{D}(X+iY)=\cal{D}X+i\cal{D}Y$. . On note $\Cun(J)$ l'ensemble
des processus $X\in \Cu (J)$ tels que pour tout
$p\in\{1,\cdots,n\}$, $\cal{D}^p X_t$ existe en tout point de $I$.
On donne \`{a} la d\'{e}finition (\ref{bonnesdiffusions}) un ensemble
$\Lambda$ qui permet de montrer que $\Cu(J)$ n'est pas trivial. En
effet $\cal{P}_{det}^1\varsubsetneq\Lambda\subset\Cu(J)$
(\textit{cf} \cite{stoc} p.26). Le calcul de $\cal{D}^p$ combine
de fa\c{c}on non triviale les quantit\'{e}s $D$ et $D_*$. \`A titre
d'exemple, on obtient sur $ \EuScript{C}^2(J)$, $\d
\cal{D}^2=\frac{DD_*+D_*D}{2}+i\frac{D^2-D_*^2}{2}$. La partie
r\'{e}elle de $\cal{D}^2$ co\"{\i}ncide donc avec
l'acc\'{e}l\'{e}ration postul\'{e}e par Nelson comme quantit\'{e}
la plus pertinente pour d\'{e}crire une notion
d'acc\'{e}l\'{e}ration pour une diffusion brownienne (cf
\cite{ne1} p.$82$).

\section{Proc\'edure de plongement stochastique}

En utilisant $\cal D$, on construit des analogues stochastiques
d'op\'{e}rateurs diff\'{e}rentiels non lin\'{e}aires.

\begin{defi}[Plongement stochastique]\label{stochastisation}
\label{stoca} On appelle plongement stochastique, relatif \`{a}
l'extension $\cal{D}_{\mu}$, d'un op\'{e}rateur $O$ qui s\'{e}crit
sous la forme : $O=a_0 (\cdot,t)+a_1 (\cdot,t)\frac{d}{dt} +\dots
+ a_n (\cdot,t)\frac{d^n}{dt^n}$ avec $a_i \in C^1(\R^d\times J)$,
$n\in\N^*$, l'op\'erateur $\cal{O} =a_0 (\cdot,t) +a_1 (\cdot,t)
{\cal D}_{\mu} +\dots + a_n (\cdot,t) {{\cal D}_{\mu}^n }$
agissant sur $\Cun(J)$.\\
Un op\'erateur $O$ \'{e}crit sous la forme : $ O=\frac{d}{dt}\circ
a(\cdot,t),$ $a\in C^1(\R^d\times J)$, est plong\'e en $
\cal{O}=\cal{D}_{\mu}\circ a(\cdot,t)$ agissant sur un
sous-ensemble de $\Cu(J)$ d\'ependant de certaines propri\'et\'es
de $a$.
\end{defi}

Un op\'erateur de la forme $ O=\frac{d}{dt}\circ a(\cdot,t)$ peut
se r\'e\'ecrire $\partial_x a(\cdot,t)\frac{d}{dt}$ qui se plonge
alors en $\partial_x a(\cdot,t)\cal{D}$. Ce dernier n'est \'egal
\`a $ \cal{O}=\cal{D}_{\mu}\circ a(\cdot,t)$ que dans certains cas
(cf \cite{stoc} p.$52$). Ceci montre en particulier que le
plongement stochastique n'est pas une application, il d\'epend du
choix d'\'ecriture de l'op\'erateur.

La notion de plongement d'op\'erateur s'\'etend de fa\c{c}on
naturelle \`a celle de plongement d'\'equation d\'efinie par un
op\'erateur $O$ d'ordre $n$ :
$O\cdot(x,\frac{dx}{dt},\cdots,\frac{d^kx}{dt^k})=0$. On d\'{e}finit
l'\'equation plong\'ee par :
$\cal{O}\cdot(X,\cal{D}_{\mu}X,\cdots,\cal{D}_{\mu}^kX)=0$ o\`u
$X\in \EuScript{C}^{n+k}(J)$. On s'occupe d\'{e}sormais du cas
lagrangien.

\begin{defi}On appelle lagrangien admissible une fonction $L: \R^d\times \C^d \to \C$ de classe $C^1$ en
sa premi\`ere variable $x$ et holomorphe en sa deuxi\`eme variable
$y$, et r\'eelle quand $y$ est r\'{e}elle. L'\'equation
\begin{equation}\label{EL}
\frac{d}{dt}\partial_yL\left(x(t),\frac{dx}{dt}(t)\right)=\partial_xL\left(x(t),\frac{dx}{dt}(t)\right)
\end{equation}
s'appelle \'equation d'Euler-Lagrange.
\end{defi}

\begin{lemm}Soit $L$ un lagrangien admissible. Le plongement stochastique de (\ref{EL}) est donn\'{e} par
\begin{equation}\label{ELS}
\cal{D}\partial_yL\left(X_t,\cal{D}X_t\right)=\partial_xL\left(X_t,\cal{D}X_t\right).
\end{equation}
\end{lemm}

On sait que l'\'equation (\ref{EL}) provient d'un principe de
moindre action (cf \cite{ar} p.$84$). Existe-il un principe de
moindre action stochastique permettant l'obtention de l'\'equation
(\ref{ELS}) ? Nous montrons dans \cite{stoc} chap.$7$ que tel est
bien le cas et on donne un lemme montrant la coh\'erence de la
proc\'{e}dure de plongement vis-\`a-vis des principes de moindre
action ainsi d\'{e}finis.

\section{\'Equation de Newton Plong\'{e}e et \'{e}quation de Schr\"{o}dinger}

Consid\'{e}rons le lagrangien admissible
$L(x,z)=\frac{1}{2}(z_1^2+\cdots+z_d^2)-U(x)$ o\`{u} $(x,z)\in
\R^d\times\C^d$ et $U$ est une fonction de classe $C^1$.
L'\'{e}quation (\ref{EL}) associ\'{e}e est l'\'{e}quation de Newton $\d
\frac{d^2x}{dt^2}(t)=-\nabla U(x(t))$. L'\'{e}quation de Newton
plong\'{e}e est alors
\begin{equation}
\cal{D}^2 X_t=-\nabla U(X_t)
\end{equation}
et co\"{\i}ncide avec l'\'{e}quation d'Euler-Lagrange
plong\'{e}e (\ref{ELS}). On se propose d'\'{e}tudier un r\'{e}sultat sur la
densit\'{e} d'un processus solution de cette \'{e}quation.

On donne dans \cite{stoc} p.$24$, suivant \cite{mns} et
\cite{thieu}, un espace sur lequel nous pourrons calculer les
d\'eriv\'ees du premier ordre $D$ et $D_*$ et les d\'eriv\'ees du
second ordre $D^2$, $DD_*$, $D_*D$ et $D_*^2$. Prenons $I=]0,1[$.
Soit $(W_t)_{t\in J}$ un mouvement brownien standard dans $\R^d$
d\'efini sur un espace probabilis\'{e} filtr\'e
$(\Omega,\cal{A},(\cal{P}_t)_{t\in J},P)$.

\begin{defi}\label{bonnesdiffusions}
On d\'esigne par $\Lambda$ l'espace des diffusions $X$
satisfaisant les conditions suivantes :
\begin{itemize}
    \item[(i)] $X$ est solution sur $J$ de l'EDS :
    $dX_t=b(t,X_t)dt+\sigma(t,X_t)dW_t,\quad X_0=X^0$
o\`u $X^0\in L^2(\Omega)$, $b:J\times \R^d\to\R^d$ et
$\sigma:J\times \R^d\to\R^d\otimes\R^d$ sont des fonctions
mesurables v\'erifiant l'hypoth\`ese : il existe une constante $K$
telle que pour tous $x,y\in\R^d$ :\\
    $\sup_t \left(\left|\sigma(t,x)-\sigma(t,y)\right|+\left|b(t,x)-b(t,y)\right|\right)\leq
    K\left|x-y\right|$ et
    $\sup_t \left(\left|\sigma(t,x)\right|+\left|b(t,x)\right|\right)\leq
    K(1+\left|x\right|)$,

    \item[(ii)] Pour tout $t\in J$, $X_t$ poss\`ede une densit\'e $p_t(x)$ en
    $x\in\R^d$,

    \item[(iii)] En posant $a_{ij}=(\sigma\sigma^*)_{ij}$, pour tout $i\in\{1,\cdots,n\}$, pour tout $t_0>0$, pour tout
    ouvert born\'e $\Xi\subset\R^d,\quad
    \int_{t_0}^1 \int_{\Xi} \left|\partial_j(a_{ij}(t,x)p_t(x))\right|dxdt <
    +\infty$,

    \item[(iv)] les fonctions $b$ et $\d (t,x)\to \frac{1}{p_t(x)}\partial_j(a_{kj}(t,x)p_t(x))$
    appartiennent \`{a} $C^1(I\times\R^d)$, sont born\'{e}es et toutes leurs d\'{e}riv\'{e}es
    du premier et second ordre sont born\'ees.
\end{itemize}
\end{defi}

On notera $\Lambda_{\sigma}$ (resp. $\Lambda^g$) le sous-ensemble
de $\Lambda$ form\'{e} par les diffusions dont le coefficient est
constant \'{e}gal \`{a} $\sigma$ (resp. dont le drift est un gradient), et
on pose $\Lambda_{\sigma}^g:=\Lambda_{\sigma}\cap\Lambda^g$.

\begin{theoreme}
Soit $X\in \Lambda$ et $f\in C^{1,2}(I\times \R^d)$ telle que
$\partial_t f$, $\nabla f$ et $\partial_{ij}f$ sont born\'{e}es. On
obtient, en adoptant la convention d'Einstein sur la sommation des
indices
\begin{eqnarray}
(\cal{D}X_t)_k & = & \left(b-\frac{1}{2p_t}\partial_j(a^{kj}p_t)+\frac{i}{2p_t}\partial_j(a^{kj}p_t)\right)(t,X_t),\\
\mathcal{D} f(t,X_t) & = & \left(\partial_t f + \mathcal{D}
X_t\cdot \nabla f +\frac{i}{2}a^{kj}\partial_{kj}f\right)(t,X_t)
.\label{deriv_fonc}
\end{eqnarray}
\end{theoreme}

On pose : $\cal{S}=\{X\in\Lambda_d \, \mid\, \cal{D}^2X_t=-\nabla
U(X(t))\}$, et pour $X\in\Lambda$ dont le drift est $b$ et la
fonction de densit\'{e} $p_t(x)$, $\Theta_X=(\R^+\times\R^d)\setminus
\{(t,x),\, \mid\, p_t(x)=0\}$.

Si $X\in \Lambda_{\sigma}^g$ alors il existe des fonctions $R$ et
$S$
diff\'{e}rentiables sur $\Theta_X$ telles que\\
$\d \cal{D}X_t=(\nabla S+i\nabla R)(X_t)$  car
$\cal{D}X_t=\left(b-\frac{\sigma^2}{2}\nabla
\log(p_t)+i\frac{\sigma^2}{2}\nabla \log(p_t)\right)(X_t)$ et $b$
est un gradient. On choisit $R(t,x)=\frac{\sigma^2}{2}\nabla
\log(p_t(x))$. Les fonctions $R$ et $S$ sont \'{e}galement
introduites par Nelson dans \cite{ne1} p.$107$.\\
On pose $A=S-iR$ et $\d \Psi_X(t,x)=e^{\frac{A(t,x)}{\sigma^2}}$.

\begin{theoreme}
\label{schro} Si $X\in\cal{S}\cap\Lambda_{\sigma}^g$, alors
$p_t(x)=|\Psi_X(t,x)|^2$ et $\Psi$ satisfait sur $\Theta_X$
l'\'{e}quation de Schr\"odinger lin\'{e}aire : $\d i\sigma^2\partial_t\Psi
+\frac{\sigma^4}{2}\Delta\Psi=U\Psi$.
\end{theoreme}

\textbf{D\'{e}monstration.} Des expressions $\d
\Psi_X(t,x)=e^{\frac{A(t,x)}{\sigma^2}}$ et
$R(t,x)=\frac{\sigma^2}{2}\nabla \log(p_t(x))$, on d\'{e}duit\\
$|\Psi_X(t,x)|^2=p_t(x)$. L'\'{e}quation de Newton plong\'{e}e peut
s'\'{e}crire $\overline{\mathcal{D}}^2X_t=-\nabla U(X_t)$ car $U$ est
r\'{e}el.\\
Or $\overline{\cal{D}}X_t=\nabla
A(t,X_t)=-i\sigma^2\frac{\nabla\Psi}{\Psi}(t,X_t)$. Donc
$-i\sigma^2\overline{\cal{D}}\frac{\nabla\Psi}{\Psi}(t,X_t)=-\nabla
U(X_t)$ et avec (\ref{deriv_fonc}) il vient\\
$i\sigma^2\left(\partial_t
\frac{\partial_k\Psi}{\Psi}+\overline{\mathcal{D}}X(t)\cdot\nabla
\frac{\partial_k\Psi}{\Psi} -i
\frac{\sigma^2}{2}\Delta\frac{\partial_k\Psi}{\Psi}\right)(t,X_t)=\partial_kU(X_t)$.
Le lemme de Schwarz donne :\\
$\overline{\cal{D}}X(t)\cdot \nabla\frac{\partial_k\Psi}{\Psi}=
-\frac{i\sigma^2}{2}\partial_k\sum_{j=1}^d\left(\frac{\partial_j\Psi}{\Psi}\right)^2$,
et $\Delta \frac{\partial_k\Psi}{\Psi}=\partial_k\sum_{j=1}^d
\left[\frac{\partial_j^2\Psi}{\Psi}-\left(\frac{\partial_j\Psi}{\Psi}\right)^2\right]$,
et par cons\'{e}quent :\\
$\ i\sigma^2\partial_k\left(\frac{\partial_t\Psi}{\Psi} -i
\frac{\sigma^2}{2}\frac{\Delta\Psi}{\Psi}\right)(t,X_t)=\partial_kU(X_t)$.
En int\'{e}grant sur $\Theta_X$ les fonctions des deux membres de la
derni\`{e}re \'{e}quation, il appara\^{\i}t des constantes qu'on peut rendre
nulles en ajoutant une fonction de $t$ convenable dans $S$. Le
r\'{e}sultat s'en d\'{e}duit.  $\square$



La partie r\'{e}elle de l'\'{e}quation de Newton plong\'{e}e co\"{\i}ncide avec
l'\'{e}quation de Newton stochastique propos\'{e}e par Nelson dans sa
th\'{e}orie dynamique des diffusions browniennes (\cite{ne1} p.$83$).
Sa partie imaginaire correspond \`{a} l'\'{e}quation $(D^2-D_*^2)X=0$.
Nous conjecturons que cette derni\`{e}re impose que le drift de $X$
doit \^{e}tre un gradient, et donc qu'il n'est pas utile de le
supposer dans le th\'{e}or\`{e}me (\ref{schro}).



\end{document}